\documentclass[11pt]{article}
\usepackage[a4paper,margin=24mm]{geometry}
\usepackage{amsmath,amssymb,amsthm,array,microtype}
\usepackage[colorlinks=true,linkcolor=blue,urlcolor=blue]{hyperref}
\newtheorem{theorem}{Theorem}

\newcommand{\blank}{\cdot}
\newcommand{\ls}[2]{\left[\begin{array}{*{#1}{c}}#2\end{array}\right]}
\newcommand{\scs}{\operatorname{scs}}

\title{A counting proof for the smallest critical set in a Latin square}
\author{Richard Bean\thanks{The mathematical argument and
writing were substantially AI-generated using OpenAI Codex and Anthropic
Claude, under my direction. Appendix~\ref{app:ai} records
these contributions, the initial prompt and the starting materials.}\\[4pt]
\small Cyber Research Centre, University of Queensland, 4072, Australia}
\date{13 October 2026}
\begin{document}
\maketitle

\begin{abstract}
A critical set in a Latin square is a set of entries that determines the
square uniquely and is minimal with this property. We prove that every
partial Latin square of order $n$ with fewer than $\lfloor n^2/4\rfloor$
filled cells has an even number of completions. Consequently every
critical set in every Latin square of order $n$ has at least
$\lfloor n^2/4\rfloor$ entries. Since the cyclic Latin square has a
critical set of exactly this size, the minimum size of a critical set
over all Latin squares of order $n$ is $\lfloor n^2/4\rfloor$. This
confirms a conjecture attributed to Nelder (1979) and made independently
by Mahmoodian and by Bate and van Rees. Hatami and Qian had shown that
this minimum is at least $n^2/10^4$ for large $n$. The proof is a
counting argument modulo two that uses only expansion of products,
cancellation in pairs and the pigeonhole principle, and it is illustrated
with partial Latin squares of orders three to seven. An appendix
discloses the substantial use of AI in developing the proof and writing
the exposition.
\end{abstract}

\section{Definitions and results}\label{sec:defs}
A \emph{Latin square} of order $n$ is an $n\times n$ array filled with
the symbols $0,1,\ldots,n-1$ so that each symbol occurs exactly once in
every row and exactly once in every column. Rows, columns and symbols
are numbered from $0$; cell $(r,c)$ lies in row $r$ and column $c$, so
the top left cell is $(0,0)$. A \emph{partial Latin square} of order $n$
is an $n\times n$ array in which some cells, the \emph{filled cells},
contain symbols, while no symbol occurs twice in a row or twice in a
column. In displays a dot marks an empty cell. A \emph{completion} of a
partial Latin square $P$ is a Latin square that agrees with $P$ in every
filled cell, and $N(P)$ is the number of completions of $P$.
Rows, columns and symbols retain their labels when completions are
counted; we do not identify different arrays under permutations.

Let $L$ be a Latin square. A \emph{defining set} of $L$ is a partial
Latin square whose filled cells all agree with $L$ and whose only
completion is $L$. A defining set is \emph{critical} if emptying any one
of its filled cells leaves a partial Latin square with more than one
completion. Let $\scs(L)$ be the smallest number of filled cells in a
critical set of $L$. This is also the smallest number of filled cells in
a defining set of $L$, because in a defining set with as few filled cells
as possible, emptying any filled cell destroys the unique completion, so
the set is critical. Let $\scs(n)$ be the minimum of
$\scs(L)$ over all Latin squares $L$ of order $n$. Finally, $B_n$ is the
\emph{cyclic square} of order $n$, whose entry in cell $(r,c)$ is
$(r+c)\bmod n$.

\paragraph{Historical context.}
The critical-set problem goes back to Nelder's 1977 newsletter
note~\cite{Nelder77}; see Bean's historical
account~\cite[Section~2.6]{BeanThesis}. The functions $\scs(n)$ and
$\operatorname{lcs}(n)$, the smallest and largest sizes of critical sets
over all Latin squares of order $n$, respectively, were introduced by
Curran and van Rees~\cite{CvR},
as explicitly credited by Stinson and van Rees~\cite[p.~442]{SvR82}.
The square-specific notation
$\scs(L)$ was introduced by Ghandehari, Hatami and
Mahmoodian~\cite{GHM}, as recorded by Bean~\cite[Section~1]{BeanTrades}.
Hatami and Qian~\cite{HQ} attribute the conjecture that
$\scs(n)=\lfloor n^2/4\rfloor$ to Nelder in a 1979 private communication
to Seberry, with later independent formulations by Mahmoodian and by
Bate and van Rees. They proved
$\scs(n)\ge 10^{-4}n^2$ for all sufficiently large
$n$~\cite[Theorem~1.2]{HQ}. The parity assertion in
Theorem~\ref{thm:parity} below implies the conjectured lower bound; it
was proposed by Bean as part~(2) of Problem BCC19.2~\cite{BeanBCC}.
The cyclic construction giving the
upper bound is classical; its sources are credited in
Section~\ref{sec:proofB}.

Related ideas occur in Balasubramanian's proof that Latin squares of
even order have an even number of transversals~\cite{Bala90}: counting
via permanent coefficients, inclusion--exclusion, equality of permanents
and determinants modulo two, and pairing complementary sets of symbols.

\begin{theorem}\label{thm:parity}
If a partial Latin square $P$ of order $n$ has fewer than
$\lfloor n^2/4\rfloor$ filled cells, then $N(P)$ is even.
\end{theorem}

\begin{theorem}\label{thm:scs}
For every $n\ge1$, \ $\scs(n)=\scs(B_n)=\lfloor n^2/4\rfloor$.
\end{theorem}

Since $1$ is odd, Theorem~\ref{thm:parity} says that a partial Latin
square with fewer than $\lfloor n^2/4\rfloor$ filled cells can never have
exactly one completion. So every defining set of every Latin square of
order $n$ has at least $\lfloor n^2/4\rfloor$ filled cells.
Theorem~\ref{thm:scs} adds that a particular choice of
$\lfloor n^2/4\rfloor$ filled cells of $B_n$ has exactly one completion.
It does not claim that $\scs(L)=\lfloor n^2/4\rfloor$ for
every individual Latin square $L$.

The proof of Theorem~\ref{thm:parity} uses only expansion of products,
cancellation in pairs, and the pigeonhole principle. We follow it
through examples of orders $3$ to $7$. The examples are taken from
cyclic squares so that they are easy to check, but the argument applies
to every partial Latin square. Sections~\ref{sec:order3}
to~\ref{sec:proofA} prove Theorem~\ref{thm:parity}, and
Section~\ref{sec:proofB} proves Theorem~\ref{thm:scs}.

\section{What we are counting: an example of order three}
\label{sec:order3}
A completion fills the dots so that every row and every column contains
each symbol exactly once. Begin with just one filled cell:
\[
 P_3=\ls{3}{0&\blank&\blank\\
                 \blank&\blank&\blank\\
                 \blank&\blank&\blank}.
\]
Here are all four completions:
\[
 \ls{3}{0&1&2\\1&2&0\\2&0&1}
 \qquad
 \ls{3}{0&1&2\\2&0&1\\1&2&0}
 \qquad
 \ls{3}{0&2&1\\1&0&2\\2&1&0}
 \qquad
 \ls{3}{0&2&1\\2&1&0\\1&0&2}.
\]
To see that the list is complete, note that the top row can be finished
in two ways. Once it is chosen, the first entry of the next row has two
possibilities, and then everything else is forced. Thus $N(P_3)=4$. The
threshold in Theorem~\ref{thm:parity} is $\lfloor 3^2/4\rfloor=2$, and
$P_3$ has one filled cell, so the theorem says that the count is even. The
argument below explains why, without listing completions.

\section{Recording the choices for each row}\label{sec:rows}
Temporarily fill the rows independently, postponing the column
conditions. For the top row of $P_3$ the two possibilities are
\[
 [\,0\ \blank\ \blank\,]
 \longrightarrow
 [\,0\ 1\ 2\,],\quad[\,0\ 2\ 1\,].
\]
Use a \emph{tag} $x_{cs}$ for the instruction ``put symbol $s$ in column
$c$''. Only placements in empty cells receive tags; filled cells do not. The
two choices above are therefore recorded by
\[
 U=x_{11}x_{22}+x_{12}x_{21}.
\]
Multiplication records choices made together, and addition lists
alternatives. The subscripts are \emph{column, symbol}, not row, column.
The nine tags form the \emph{tag table} $X$, whose rows are the columns
$c$ of the Latin square and whose columns are the symbols $s$. It does not
depend on $P_3$: it lists all nine column-symbol pairs, without specifying
a Latin-square row. It is a bookkeeping table, not a Latin square. The
partial square only decides which part of $X$ each row may use. Here $X$ is shown
next to $P_3$ for comparison:
\[
 \underbrace{\ls{3}{0&\blank&\blank\\
                 \blank&\blank&\blank\\
                 \blank&\blank&\blank}}_{\text{Latin-square rows and columns}}
 \qquad
 X=\begin{array}{c|ccc}
      &s=0&s=1&s=2\\\hline
 c=0&x_{00}&x_{01}&x_{02}\\
 c=1&x_{10}&x_{11}&x_{12}\\
 c=2&x_{20}&x_{21}&x_{22}
 \end{array}.
\]
Each of the other two rows of $P_3$ is empty, so each has six possible
fillings, all recorded by the same expression
\begin{align*}
 H={}&x_{00}x_{11}x_{22}+x_{00}x_{12}x_{21}
       +x_{01}x_{10}x_{22}\\
    &+x_{01}x_{12}x_{20}+x_{02}x_{10}x_{21}
       +x_{02}x_{11}x_{20}.
\end{align*}
The product $UHH$ lists all $2\cdot6\cdot6=72$ ways of finishing the
rows independently. Most of them fail the column conditions. For
example, this filling has valid rows but is \emph{not} a Latin square:
\[
 \ls{3}{0&1&2\\1&0&2\\1&2&0}.
\]
The two new $1$s in column $0$ give two copies of $x_{01}$ in its
product of tags: a repeated tag detects a repeated symbol in a column.

The filled cell already supplies the pair ``column $0$, symbol $0$''. A
completion must supply every other column-symbol pair exactly once, so
its eight new placements give exactly
\[
 V=x_{01}x_{02}x_{10}x_{11}x_{12}x_{20}x_{21}x_{22}.
\]
An \emph{expression} is a sum of products of tags. A product may have one
tag or several, a product of no tags is the number $1$, and the same
product may occur more than once; the number of times it occurs is its
\emph{coefficient}. The order of tags in a product does not matter:
$x_{11}x_{22}$ and $x_{22}x_{11}$ are the same product.
For example, each term of $U$ is a product of two
tags, and $1+x_{00}$ is the number $1$ plus a single tag. When $UHH$ is
multiplied out, $V$ occurs once for each completion, so its coefficient
is $4$.

\paragraph{Matching sums and row expressions in general.}
Let $Y$ be a square table of expressions, such as the tag table $X$.
Choose a set $F$ of its rows and a set $G$ of its columns with the same
number $t$ of elements, and consider every way of choosing one entry in
each row of $F$ so that each column of $G$ is used exactly once. The
\emph{matching sum} $M_Y(F,G)$ is the sum, over all these ways, of the
product of the chosen entries. Its \emph{size} is $t$, and a matching
sum of size $0$ is the empty product $1$. For example, the top row of
$P_3$ has empty columns $1$ and $2$ and is missing symbols $1$ and $2$,
so its choices are the matching sum on rows $\{1,2\}$ and columns
$\{1,2\}$ of $X$: $U=M_X(\{1,2\},\{1,2\})$. An empty row may put any
symbol in any column, so its choices are the matching sum of the whole
of $X$: $H=M_X(\{0,1,2\},\{0,1,2\})$.

Now let $P$ be any partial Latin square of order $n$ with $m$ filled
cells. For each row $r$, let $d_r$ be its number of filled cells, $F_r$
the set of its
empty columns, and $G_r$ the set of symbols missing from it; both sets
have $t_r=n-d_r$ elements. The \emph{row expression} $M_X(F_r,G_r)$
lists all ways to put the missing symbols into the empty cells of row
$r$, so the product $\prod_rM_X(F_r,G_r)$ of the row expressions lists
all ways to fill every row independently. Call $(c,s)$ a \emph{supplied pair}
if some filled cell puts symbol $s$
in column $c$; different filled cells give different supplied pairs,
since no symbol is repeated in a column. Let $V$ be the product of the
tags $x_{cs}$ of all pairs $(c,s)$ that are not supplied pairs. It has
$n^2-m$ distinct tags, the same number as the empty cells. A product of
row choices equals
$V$ exactly when every pair that is not supplied is used once and no
supplied pair is used again, that is, exactly when every column contains
each symbol once. Such a choice of rows is precisely a completion, and
different completions are different choices. Therefore
\begin{equation}\label{eq:coef}
 \text{the coefficient of $V$ in the ordinary expansion of }
 \prod_r M_X(F_r,G_r)\ \text{is } N(P).
\end{equation}
A product in this expansion that repeats no tag but differs from $V$ is
not counted: it
contains the tag of a supplied pair, so some column would contain the
symbol of a filled cell twice.

\section{Two rules for keeping only the parity}\label{sec:rules}
The coefficient of $V$ is the number of completions $N(P)$. If this
number is even, it is either zero or at least two, so $P$ cannot have
exactly one completion. We therefore need only its parity, not its exact
value, and use two bookkeeping rules:
\begin{enumerate}
\item Cancel two identical products whenever they occur. This reduces
that product's coefficient by two, preserving its parity. In particular,
two occurrences of $V$ record two completions and may be cancelled
together.
\item Discard any product that contains some tag twice. In the original
row-choice expansion, a repeated tag $x_{cs}$ means that symbol $s$ has
been placed twice in column $c$, so that filling is not a Latin square.
Such a product cannot be the target $V$, which contains no repeated tag.
\end{enumerate}
Both rules therefore preserve the coefficient of $V$ modulo $2$.
The rules can be applied at any stage of a calculation
without changing the final result: later multiplication cannot remove a
repeated tag, and two equal terms stay equal after the same further
operations. Write $\equiv$ for equality after applying the rules. These
are rules for a tally, not statements about numerical values of the
tags; for example, $x_{11}^2\equiv0$ because that product is discarded.

An expression is \emph{constant-free} if every one of its terms contains at
least one tag, that is, if it has no constant term. For example,
squaring the constant-free expression $U$ gives
\[
 U^2=(x_{11}x_{22})^2+(x_{12}x_{21})^2
                    +2x_{11}x_{22}x_{12}x_{21}\equiv0:
\]
the first two terms repeat tags, and the two copies of the cross term
cancel. The same reasoning applies to every constant-free expression $a$: the
square of each term of $a$ repeats a tag, and each product of two
different terms occurs twice. Since $(1+a)^2=1+2a+a^2$, it follows that
\[
 \boxed{\begin{gathered}
 a^2\equiv0\quad\text{and}\quad(1+a)^2\equiv1\\
 \text{for every constant-free expression }a.
 \end{gathered}}
\]
In particular, a product in which the same constant-free expression occurs
twice is $\equiv0$. This holds even when $a$ is a long sum of products
of tags, which matters when the entries of a table become more
complicated.

\section{Adding dummy entries without changing the count}
\label{sec:dummy}
Keep the same partial Latin square, but change one entry of its tag
table:
\[
 P_3=\ls{3}{0&\blank&\blank\\
                 \blank&\blank&\blank\\
                 \blank&\blank&\blank},
 \qquad
 Y=\ls{3}{1+x_{00}&x_{01}&x_{02}\\
                x_{10}&x_{11}&x_{12}\\
                x_{20}&x_{21}&x_{22}}.
\]
The added $1$ is a bookkeeping option, called a \emph{dummy}~$1$. It is
not a new filled cell, and it is not the Latin-square symbol $1$. Form the row
expressions from the changed table $Y$ by the same rule as before, and
only then apply the bookkeeping rules. The tags themselves and the two
rules remain unchanged; only entries of the table are altered.

Each term of the product of the three row expressions still multiplies
eight entries, one for each empty cell. If it takes the dummy $1$ from
any of them, the resulting product has at most seven tags and cannot be
$V$, which has eight. So adding the $1$ does not change the coefficient
of $V$. It does change other coefficients, which we do not need.

In general, add $1$ to the $k=\lfloor n/2\rfloor$ entries of the tag
table on the main diagonal in positions $(j,j)$ for $0\le j<k$. These positions lie in
distinct rows and distinct columns. We call these entries
\emph{marked}: each marked entry is $1$ plus a constant-free expression, and
every other entry is constant-free. Each term of the product of the row
expressions multiplies $n^2-m$ entries, and $V$ has $n^2-m$ tags, so
taking a dummy $1$ anywhere leaves too few tags. Hence the coefficient
of $V$ in the product of the row expressions formed from $Y$ is still
$N(P)$, and after the rules it is $N(P)$ modulo $2$.

\section{Removing a dummy entry}\label{sec:remove}

\subsection{The calculation at order three}\label{sec:remove3}
Put $p=1+x_{00}$, so that $p^2\equiv1$. The following operations change
the tag table $Y$, not the partial Latin square.

Row $0$ of $Y$ is $(p,\,x_{01},\,x_{02})$. Add $x_{10}p$ times row $0$
to row $1$, and $x_{20}p$ times row $0$ to row $2$. Row $1$ becomes
\[
 \bigl(x_{10}+x_{10}p\,p,\ \ x_{11}+x_{10}p\,x_{01},\ \
       x_{12}+x_{10}p\,x_{02}\bigr),
\]
and its first entry is $x_{10}+x_{10}p\,p\equiv x_{10}+x_{10}\equiv0$,
because $p\,p\equiv1$. Row $2$ changes in the same way, with $x_{20}$ in
place of $x_{10}$. So the first column becomes $(p,0,0)$.

Next add $px_{01}$ times column $0$ to column $1$, and $px_{02}$ times
column $0$ to column $2$. Column $0$ is now $(p,0,0)$, so these additions
change only row $0$, where $x_{01}+p\,x_{01}\,p\equiv0$ and
$x_{02}+p\,x_{02}\,p\equiv0$. Altogether,
\[
 \ls{3}{p&x_{01}&x_{02}\\x_{10}&x_{11}&x_{12}\\x_{20}&x_{21}&x_{22}}
 \quad\longrightarrow\quad
 \ls{3}{p&0&0\\0&z_{11}&z_{12}\\0&z_{21}&z_{22}}.
\]
Row $0$ and column $0$ now contain only $p$. Deleting them leaves the
\emph{remaining table} $Z$, the $2\times2$ table in rows $1,2$ and
columns $1,2$. We keep these labels, so $z_{vs}$ denotes the entry of $Z$
in row $v$ and column $s$. By the row additions above,
\[
 z_{vs}=x_{vs}+x_{v0}\,p\,x_{0s}\qquad(v,s\in\{1,2\}),
\]
so that
\[
 Z=\ls{2}{z_{11}&z_{12}\\z_{21}&z_{22}}
  =\ls{2}{x_{11}+x_{10}px_{01}&x_{12}+x_{10}px_{02}\\
             x_{21}+x_{20}px_{01}&x_{22}+x_{20}px_{02}}.
\]
To see the tags in an entry of $Z$, substitute $p=1+x_{00}$ and multiply
out. For example,
\[
 z_{11}=x_{11}+x_{10}(1+x_{00})x_{01}
       =x_{11}+x_{10}x_{01}+x_{10}x_{00}x_{01}.
\]
Each of these three terms contains at least one tag, so $z_{11}$ has no
constant term: it is constant-free. The same holds for every entry of $Z$, since
$z_{vs}=x_{vs}+x_{v0}x_{0s}+x_{v0}x_{00}x_{0s}$. So using the same entry
of $Z$ twice gives $0$ under the rules. Section~\ref{sec:three-rows}
rewrites the row expressions of $P_3$ in terms of the entries of $Z$.
After multiplying these rewritten expressions, expand the matching sums
while keeping each entry of $Z$ intact. If every resulting product uses
some entry of $Z$ at least twice, the whole product is $\equiv0$.
Because the rules preserve the parity of the coefficient of $V$, the
number of completions $N(P_3)$ is then even. From now on, the size of a
matching sum of $Z$ counts entries taken from $Z$, each of which may
contain several tags.

\subsection{Row and column additions}\label{sec:additions}
We now show that these additions change a matching sum only by adding
multiples of other matching sums of the same size. In the smallest case, a
two-row matching sum is $ad+bc$, and adding $e$ times the first row to
the second gives
\[
 \ls{2}{a&b\\c&d}
 \longrightarrow
 \ls{2}{a&b\\c+ea&d+eb},
 \qquad
 a(d+eb)+b(c+ea)=ad+bc+2eab\equiv ad+bc.
\]
In general, suppose $e$ times row $u$ of a table is added to a different
row $v$, and consider a matching sum of the new table that selects row
$v$. Each of its terms takes one entry $y_{vs}+e\,y_{us}$ from row $v$,
so the new matching sum is the old one plus an \emph{extra part}, in
which the entry from row $v$ is replaced by $e$ times the entry of row
$u$ in the same column.
\begin{itemize}
\item If row $u$ is not selected, the extra part is $e$ times the
matching sum that selects row $u$ instead of row $v$, with the same
columns. It has the same size.
\item If row $u$ is also selected, each extra term contains
$e\,y_{us}\,y_{us'}$, where $s$ is the column used for row $v$ and $s'$
the column used for row $u$. Exchanging $s$ and $s'$ gives an equal
term, so the extra terms cancel in pairs.
\end{itemize}
A matching sum that does not select row $v$ does not change. So every
matching sum of the new table is $\equiv$ a sum of expressions times
matching sums of the old table, of the same size. Adding $e$ times row
$u$ to row $v$ a second time restores the old table, because
$2e\,y_{us}\equiv0$. Here $e$ is held fixed as an expression in the
original tags, even if it was calculated from entries of the table.
Applying the same reasoning to that second
addition, every matching sum of the old table is $\equiv$ a sum of
expressions times matching sums of the new table, of the same size.
Column additions behave in the same way, with rows and columns
exchanged.

\subsection{The general removal step}\label{sec:general-removal}
The same two stages remove a marked entry from any table. Let $Y$ be a
square table whose entry $p=y_{ab}$, in row $a$ and column $b$, is $1$
plus a constant-free expression, so that $p\,p\equiv1$.
\begin{enumerate}
\item For every row $v\ne a$, add $y_{vb}\,p$ times row $a$ to row $v$.
The entry of row $v$ in column $b$ becomes
$y_{vb}+y_{vb}\,p\,p\equiv0$, and each other entry of row $v$ becomes
\[
 z_{vs}=y_{vs}+y_{vb}\,p\,y_{as}\qquad(v\ne a,\ s\ne b).
\]
\item For every column $s\ne b$, add $p\,y_{as}$ times column $b$ to
column $s$. Column $b$ is now zero except for $p$, so this changes only
row $a$, whose entry in column $s$ becomes $y_{as}+p\,y_{as}\,p\equiv0$.
\end{enumerate}
The result has $p$ in position $(a,b)$, zeros in the rest of row $a$
and column $b$, and the smaller table $Z=(z_{vs})$ in the remaining rows
and columns. A matching sum of the result that selects both row $a$ and
column $b$ must take $p$ from row $a$, so it is $p$ times a matching sum
of $Z$ of size one less. One that selects neither is a matching sum of
$Z$ of the same size. One that selects exactly one of them includes a
row or column of zeros, so it is $\equiv0$. Together with
Section~\ref{sec:additions}, this proves:
\begin{quote}
\emph{every matching sum of $Y$ is $\equiv$ a sum of expressions times
matching sums of $Z$, each of the same size or one less.}
\end{quote}
All the multiplying expressions are still expressions in the original
tags. The entries of $Z$ are also expressions in those tags, not new tags.
We say that the marked entry has been \emph{removed}. Now suppose, as in
all our tables, that the other marked entries of $Y$ are $1$ plus constant-free
expressions, lying in rows other than $a$ and columns other than $b$,
and that every unmarked entry is constant-free. For $v\ne a$ the entry $y_{vb}$
is unmarked, because no other marked entry lies in column $b$. So it is
constant-free, and each correction $y_{vb}\,p\,y_{as}$ is constant-free. Hence the
remaining marked entries of $Z$ are still $1$ plus constant-free expressions,
every other entry of $Z$ is still constant-free, and the next marked entry can
be removed in the same way.

\subsection{What happens to the three row expressions}
\label{sec:three-rows}
Return to $P_3$, where one marked entry has been removed, so that each
size can drop by at most one:
\[
 \ls{3}{0&\blank&\blank\\
                 \blank&\blank&\blank\\
                 \blank&\blank&\blank}
 \qquad\begin{array}{c|ccc}
       \text{Latin-square row }r&0&1&2\\\hline
       t_r&2&3&3\\
       \text{size left, at least}&1&2&2
       \end{array}.
\]
For an explicit check, put $D=z_{11}z_{22}+z_{12}z_{21}$. Each empty row
of $P_3$ has the full matching sum of $Y$ as its row expression. Every
row and column is selected, so by Section~\ref{sec:additions} the
additions preserve it under the bookkeeping rules. Thus it is
$\equiv pD$, the full matching sum of the transformed $3\times3$ table.
The top row's expression satisfies
\begin{align*}
 U\equiv D+p(&x_{10}x_{01}z_{22}+x_{10}x_{02}z_{21}\\
             &+x_{20}x_{01}z_{12}+x_{20}x_{02}z_{11}),
\end{align*}
which can be checked by substituting the definitions of the $z$'s and
cancelling equal pairs. Its first part uses two entries of $Z$, and
every other part uses one, times an expression. So the sizes $1$, $2$,
$2$ are visible in the formulas.

After multiplying the three transformed row expressions, each term uses
at least $1+2+2=5$ entries of $Z$, but $Z$ has only four positions. Some
entry is used twice, and its square is $0$ under the rules; other
factors cannot change this. So the whole product is $\equiv0$, the
coefficient of $V$ is even, and indeed $N(P_3)=4$.

\section{A balanced example of order five}\label{sec:order5}
Consider five filled cells, one in each row and column:
\[
 P_5=\ls{5}{0&\blank&\blank&\blank&\blank\\
 \blank&2&\blank&\blank&\blank\\
 \blank&\blank&4&\blank&\blank\\
 \blank&\blank&\blank&1&\blank\\
 \blank&\blank&\blank&\blank&3}.
\]
These filled cells agree with the cyclic square $B_5$. Their number,
five, is below the threshold $\lfloor 25/4\rfloor=6$, and every row has
four empty cells. Add dummy $1$s at the tag-table positions $(0,0)$ and $(1,1)$, and
remove them one after the other. The second removal is allowed by the
last paragraph of Section~\ref{sec:general-removal}: the two marked
positions are in different rows and columns, so the first removal
changes the other entries only by constant-free corrections, and the entry at
$(1,1)$ is still $1$ plus a constant-free expression.

After both removals only a $3\times3$ table remains:
\[
 \underbrace{\ls{3}{z_{22}&z_{23}&z_{24}\\
                         z_{32}&z_{33}&z_{34}\\
                         z_{42}&z_{43}&z_{44}}}_{9\text{ constant-free entries}}.
\]
Each of the five row expressions still uses at least $4-2=2$ entries of
this \emph{same} table, so their product uses at least $5\cdot2=10$
entries from nine positions. A repeated entry forces cancellation, and
the number of completions is even. (Direct enumeration gives $48$; the
argument needs only evenness.) This example has no particularly sparse
row: the argument uses the total number of entries required by all the
rows together.

\section{Orders four, six and seven}\label{sec:more}
At order $4$, take these three filled cells:
\[
 P_4=\ls{4}{0&\blank&\blank&\blank\\
 \blank&2&\blank&\blank\\
 \blank&\blank&0&\blank\\
 \blank&\blank&\blank&\blank}.
\]
The row expressions have sizes $(3,3,3,4)$. Two removals leave sizes at
least $(1,1,1,2)$, so five entries must come from a $2\times2$ table, and
the completion count is even. Direct enumeration gives $8$.

At order $6$, the threshold is $9$. Here is an example with $8$ filled cells:
\[
 P_6=\ls{6}{0&\blank&\blank&3&\blank&\blank\\
 \blank&2&\blank&\blank&5&\blank\\
 \blank&\blank&4&\blank&\blank&\blank\\
 \blank&\blank&\blank&0&\blank&\blank\\
 \blank&\blank&\blank&\blank&2&\blank\\
 \blank&\blank&\blank&\blank&\blank&4}.
\]
The row sizes are $(4,4,5,5,5,5)$. Three removals leave at least
$(1,1,2,2,2,2)$, whose sum is $10$. The remaining table has nine
positions, so again the product vanishes under the rules. (Direct
enumeration gives $812$ completions.)

At order $7$, take $11$ filled cells, just below the threshold $12$:
\[
 P_7=\ls{7}{0&\blank&\blank&3&\blank&\blank&\blank\\
 \blank&2&\blank&\blank&5&\blank&\blank\\
 \blank&\blank&4&\blank&\blank&0&\blank\\
 \blank&\blank&\blank&6&\blank&\blank&2\\
 \blank&\blank&\blank&\blank&1&\blank&\blank\\
 \blank&\blank&\blank&\blank&\blank&3&\blank\\
 \blank&\blank&\blank&\blank&\blank&\blank&5}.
\]
Every row and every column has one or two filled cells, and every
symbol occurs in one or two of them.
The row sizes are $(5,5,5,5,6,6,6)$. After three removals they are at
least $(2,2,2,2,3,3,3)$, totalling $17$, while the remaining $4\times4$
table has only $16$ positions. The balanced placement of the filled cells does
not prevent the repetition. (Direct enumeration gives $47{,}272$
completions.)

The comparisons for the examples so far are:
\[
\begin{array}{c|c|c|c|c}
 n&\text{filled cells}&\text{removals }k&
       \text{entries required, at least}&\text{positions available}\\\hline
 3&1&1&5&4\\
 4&3&2&5&4\\
 5&5&2&10&9\\
 6&8&3&10&9\\
 7&11&3&17&16
\end{array}
\]
Every displayed partial square is contained in its cyclic square, so it
has at least one completion. Since its number of completions is even,
it has at least two.

\section{Filled cells concentrated in a few rows}\label{sec:concentrated}
At order $5$, move from one filled cell in each row to a completely
filled first row, keeping five filled cells:
\[
 Q_5=\ls{5}{0&1&2&3&4\\
 \blank&\blank&\blank&\blank&\blank\\
 \blank&\blank&\blank&\blank&\blank\\
 \blank&\blank&\blank&\blank&\blank\\
 \blank&\blank&\blank&\blank&\blank}.
\]
Now the sizes are $(0,5,5,5,5)$. After two removals the lower bounds are
$(0,3,3,3,3)$, totalling $12$, again more than nine positions. A matching
sum of size $0$ is the empty product $1$: a completed row has one way to
do nothing, and its size never becomes negative. This is why the lower
bound after $k$ removals is $\max(t_r-k,0)$ rather than $t_r-k$, and the
same rule covers balanced and concentrated placements. (This square has
$1344$ completions.)

\section{Proof of Theorem A for every order}\label{sec:proofA}
Let $P$ be a partial Latin square of order $n$ with
$m<\lfloor n^2/4\rfloor$ filled cells. For $n=1$ the threshold is $0$, so there
is nothing to prove; let $n\ge2$ and $k=\lfloor n/2\rfloor$. Then
$k(n-k)=\lfloor n^2/4\rfloor$: if $n=2k$ both sides equal $k^2$, and if
$n=2k+1$ both equal $k^2+k$.

Add the $k$ dummy $1$s of Section~\ref{sec:dummy}, so that the
coefficient of $V$ in the product of the row expressions is still
$N(P)$. Remove the marked entries one at a time as in
Section~\ref{sec:general-removal}. Each removal turns every matching sum
into a sum of expressions times matching sums of the smaller table, of
the same size or one less, and leaves the remaining marked entries of
the form $1$ plus a constant-free expression. After $k$ removals, the product of
the row expressions is $\equiv$ a sum of terms, each an expression times
a product of matching sums of one final $(n-k)\times(n-k)$ table of
constant-free entries. In each such product the factor coming from row $r$ has
size at least $\max(t_r-k,0)$, since a size drops by at most one per
removal and is never negative. Expand these matching sums, keeping each
entry of the final table as a single expression. Every resulting product
uses at least
\[
 W=\sum_r\max(t_r-k,0)\ \ge\ \sum_r(t_r-k)\ =\ n(n-k)-m
\]
entries of the final table. Since $m<k(n-k)$,
\[
 W>n(n-k)-k(n-k)=(n-k)^2,
\]
which is the number of positions in the final table. By the pigeonhole
principle, every term uses some entry $z$ at least twice and therefore
contains the factor $z^2\equiv0$, since $z$ is constant-free. This remains true
after multiplication by any other expressions.
Hence the product of the row expressions is $\equiv0$. The bookkeeping
rules preserve the parity of the coefficient of $V$, whose value before
applying the rules is $N(P)$. Therefore $N(P)$ is even. This proves
Theorem~\ref{thm:parity}.

For a partial Latin square with no completion this says only that $0$
is even. If the filled cells are taken from a Latin square, there is already one
completion, so evenness forces a second.

\section{At the threshold, and proof of Theorem B}\label{sec:proofB}
The strict inequality in Theorem~\ref{thm:parity} matters. Fill a second cell in the example of order three:
\[
 C_3=\ls{3}{0&\blank&\blank\\
                 \blank&\blank&\blank\\
                 \blank&\blank&1}
 \quad\longrightarrow\quad
 \ls{3}{0&1&2\\1&2&0\\2&0&1}=B_3.
\]
Only the first of the four completions of $P_3$ has symbol $1$ in cell $(2,2)$,
so $N(C_3)=1$. Here the sizes are $(2,3,2)$, and one removal leaves at
least $(1,2,1)$, totalling exactly four. Four entries from four
positions need not repeat a position, so the counting argument correctly
stops at the threshold.

\subsection*{The corner sets}
This is the standard cyclic-square construction, written with indices
starting at zero. Its history includes Curran and van Rees~\cite{CvR}
and Smetaniuk~\cite{Smetaniuk}; Cooper, Donovan and Seberry~\cite{CDS}
established the required criticality in odd order. Bate and van
Rees~\cite[p.~49]{BvR03} give the construction and this attribution.
The construction is included here to make the upper-bound argument
self-contained.

For any $n$, let $k=\lfloor n/2\rfloor$ and keep the cells $(r,c)$ of
$B_n$ whose ordinary integer sum satisfies
\[
 r+c<k\qquad\text{or}\qquad r+c\ge n+k.
\]
Call the resulting partial Latin square $C_n$. It consists of two corner
triangles, and for $n=3$ it is the square $C_3$ above. For $n=5$ the
conditions are $r+c<2$ or $r+c\ge7$:
\[
 C_5=\ls{5}{0&1&\blank&\blank&\blank\\
 1&\blank&\blank&\blank&\blank\\
 \blank&\blank&\blank&\blank&\blank\\
 \blank&\blank&\blank&\blank&2\\
 \blank&\blank&\blank&2&3}.
\]
There are three filled cells in each corner, six altogether. Fill the empty
cells row by row from top to bottom, and from right to left within each
row. The top right cell $(0,4)$ cannot be $0$ or $1$, which are already
in its row, or $2$ or $3$, which are already in its column. It must be
$4$:
\[
 \ls{5}{0&1&\blank&\blank&\boxed{4}\\
 1&\blank&\blank&\blank&\blank\\
 \blank&\blank&\blank&\blank&\blank\\
 \blank&\blank&\blank&\blank&2\\
 \blank&\blank&\blank&2&3}.
\]
At $(0,3)$, the row excludes $0,1,4$ and the column excludes $2$,
forcing $3$. Then the last empty cell of the top row must be $2$:
\[
 \ls{5}{0&1&\boxed{2}&\boxed{3}&4\\
 1&\blank&\blank&\blank&\blank\\
 \blank&\blank&\blank&\blank&\blank\\
 \blank&\blank&\blank&\blank&2\\
 \blank&\blank&\blank&2&3}.
\]
In the next row, cell $(1,4)$ is forced to be $0$: its row excludes $1$
and its column excludes $2,3,4$. Moving left then forces $4,3,2$, in
that order:
\[
 \ls{5}{0&1&2&3&4\\
 1&\boxed{2}&\boxed{3}&\boxed{4}&\boxed{0}\\
 \blank&\blank&\blank&\blank&\blank\\
 \blank&\blank&\blank&\blank&2\\
 \blank&\blank&\blank&2&3}.
\]
Continuing in the same order forces the remaining entries and gives
\[
 B_5=\ls{5}{0&1&2&3&4\\1&2&3&4&0\\2&3&4&0&1\\
                  3&4&0&1&2\\4&0&1&2&3}.
\]

\subsection*{Why the same order works for every $n$}
Fill the empty cells of $C_n$ row by row from the top, and from right to
left within each row. Suppose that every cell filled so far agrees with
$B_n$, and consider the next empty cell $(r,c)$, so that
$k\le r+c<n+k$. Its intended symbol $(r+c)\bmod n$ is not yet present in
its row or column, because everything present so far agrees with $B_n$.
Take any other symbol. In $B_n$ it lies in row $r$ at some column
$j\ne c$.
\begin{itemize}
\item If $j>c$, then cell $(r,j)$ has already been forced or is a filled
cell of $C_n$, so the symbol is already in row $r$.
\item If $j<c$ and $(r,j)$ is a filled cell of $C_n$, the symbol is again
already in row $r$.
\item Otherwise $j<c$ and $(r,j)$ is not a filled cell of $C_n$, so
$u=r+j$ satisfies
$k\le u<r+c<n+k$. In $B_n$ the symbol lies in column $c$ at row
$i=(u-c)\bmod n$. If $u-c\ge0$, then $i=r+j-c<r$, so cell $(i,c)$ lies
in an earlier row and has already been forced or is a filled cell of
$C_n$. If $u-c<0$, then $i=u-c+n$ and $i+c=u+n\ge n+k$, so $(i,c)$ is a
filled cell of $C_n$. Either
way the symbol is already in column $c$.
\end{itemize}
So every other symbol is excluded and the cell is forced to its
intended symbol. By induction every empty cell is forced, and $C_n$ is
a defining set of $B_n$. (For $n=1$, $C_1$ has no filled cells and the single
cell is forced trivially.)

\subsection*{The number of filled cells}
For $0\le u<k$, the cells with $r+c=u$ are
$(0,u),(1,u-1),\ldots,(u,0)$, so there are $u+1$ of them, and these
diagonals contribute $1+2+\cdots+k=k(k+1)/2$ filled cells. For
$n+k\le u\le 2n-2$, the cells with $r+c=u$ have
$u-n+1\le r\le n-1$, so there are $2n-1-u$ of them, and these diagonals
contribute $1+2+\cdots+(n-k-1)=(n-k-1)(n-k)/2$ filled cells. If $n=2k$, the
total is $k(k+1)/2+(k-1)k/2=k^2$; if $n=2k+1$, it is
$k(k+1)/2+k(k+1)/2=k(k+1)$. In both cases $C_n$ has
$k(n-k)=\lfloor n^2/4\rfloor$ filled cells. For orders $3,4,5,6,7$ this
gives $2,4,6,9,12$ filled cells.

\subsection*{Proof of Theorem B}
For $n=1$ the empty set has exactly one completion and is critical, so
the result holds. Suppose $n\ge2$.
A defining set has exactly one completion, and $1$ is odd, so by
Theorem~\ref{thm:parity} it has at least $\lfloor n^2/4\rfloor$ filled
cells.
Hence $\scs(L)\ge\lfloor n^2/4\rfloor$ for every Latin square $L$ of
order $n$, and $\scs(n)\ge\lfloor n^2/4\rfloor$. The corner set $C_n$ is
a defining set of $B_n$ with exactly $\lfloor n^2/4\rfloor$ filled cells.
It is critical: emptying any filled cell leaves $\lfloor n^2/4\rfloor-1$
filled cells, so by
Theorem~\ref{thm:parity} an even number of completions, at least one of
which is $B_n$; hence at least two. Therefore
$\scs(B_n)\le\lfloor n^2/4\rfloor$, and
\[
 \lfloor n^2/4\rfloor\le\scs(n)\le\scs(B_n)\le\lfloor n^2/4\rfloor.
\]
Thus
\[
 \boxed{\scs(n)=\scs(B_n)=\lfloor n^2/4\rfloor.}
\]
Here $\scs(n)$ is the minimum over all Latin squares of order $n$; the
argument does not classify all Latin squares or all critical sets
attaining this minimum.

\appendix
\section{AI contribution and development record}\label{app:ai}
\subsection{Contributions and checking}
The proof and the writing in this document were substantially generated
by AI. During iterative work on 11--13 October 2026, I chose
the problem, supplied background papers, directed the investigation
through prompts and feedback, and requested reviews, corrections and
simpler explanations. OpenAI Codex and Anthropic Claude were used to
generate and develop mathematical arguments, draft and revise the
exposition, and review and refine successive versions. Their contribution
extended to the substance of the proof, not just its wording or
typesetting. The Codex session records identify its model as
\texttt{gpt-6-astra}; all Claude sessions used Claude Opus 5.5 (model
identifier \texttt{claude-opus-5-5}).

Codex also generated and ran programs checking small examples and
intermediate identities. Claude independently ran exhaustive completion
counts at orders $3$, $4$ and $5$ (every partial Latin square with
$\lfloor n^2/4\rfloor-1$ filled cells), randomized tests at orders $6$ and $7$,
randomized tests of the identities in Section~\ref{sec:remove}, and
recounts of every example.

\subsection{The initial prompt}
The initial prompt was sent on 11 October 2026.
The prefix \texttt{/goal} was an interface
command.
\begin{quote}
\small\ttfamily\raggedright
/goal prove that the smallest critical set of a Latin square of order n
is of size floor(n\textasciicircum{}2/4).
\medskip
\noindent
the papers in the directory maths will be good background. if you cannot
prove that then at least try to prove if for just the special case of
Z\_n with n odd; or improve the best known bound.
\end{quote}

Here \texttt{Z\_n} denotes the cyclic square called $B_n$ in this
paper. The target minimum is $\scs(n)$, as defined in
Section~\ref{sec:defs}.

\subsection{Background supplied at the outset}
With the initial prompt, I supplied two papers by Bate and
van Rees~\cite{BvR99,BvR03}, and papers by Cavenagh~\cite{Cavenagh},
Ghandehari, Hatami and Mahmoodian~\cite{GHM}, and Hatami and
Qian~\cite{HQ}. The roles of these sources are described in
Appendix~\ref{app:sources}.

\subsection{Subsequent development}
The initial prompt was followed by many rounds of investigation,
feedback and revision. Further prompts asked for trade and cycle-switch
arguments in odd cyclic squares, for constraints on sparse and balanced
sets of filled cells, and for a trade-covering integer programme for order
11. The investigation then returned to the general order-$n$ problem
and the stronger parity statement of Theorem~\ref{thm:parity}. An algebraic version of
the argument, produced in the Codex sessions, preceded the present proof.
It relied on an external theorem on determinantal ideals. I
requested checks and corrections using Codex and Claude, and then asked
for a version without the algebraic machinery. In response, Claude found
the self-contained step used here (the dummy entries of
Section~\ref{sec:dummy} and their removal in Section~\ref{sec:remove}),
which replaces that external theorem. The exposition was then revised by
Codex and Claude, and explicit small Latin-square examples were added at
my request. The initial prompt was not the sole input to the
mathematical text above.

Appendix~\ref{app:sources} records the
roles of the sources and distinguishes direct consultation from historical
attribution. This appendix gives an excerpt from the initial prompt and a summary of
the subsequent process, rather than a verbatim transcript of all the
Codex and Claude exchanges.

\section{Sources and their roles in the development}\label{app:sources}
The following account distinguishes sources consulted during the
investigation from historical references added when preparing this
version. It records relevant uses rather than treating every search
result or briefly inspected paper as a source of the proof.

\paragraph{Latin-square background actually consulted.}
Hatami and Qian~\cite{HQ} supplied the earlier quadratic lower-bound
starting point. The two papers of Bate and van Rees~\cite{BvR99,BvR03}
supplied the strong-critical-set result, cyclic constructions and
small-order background. Cavenagh~\cite{Cavenagh} and
Bean~\cite{BeanTrades} were used in the earlier investigations of
sparse rows, trades and trade-covering computations.
Ghandehari, Hatami and Mahmoodian~\cite{GHM} were consulted for the
distinction between $\scs(n)$ and the behaviour of $\scs(L)$ across
individual squares. Bean's problem~\cite{BeanBCC} supplied the precise
parity assertion targeted here.

\paragraph{Sources used for the earlier algebraic argument.}
Bruns and Vetter~\cite{BrunsVetter}, particularly Theorem~2.10 and
Remark~2.12, supplied an external theorem used in the earlier proof.
The proof of Theorem~1.2 in Miller, Singh and Varbaro~\cite{MSV} was
consulted for the related method, and Hochster and
Huneke~\cite{HH} were consulted during review of that version.
These works belong in the development record even though the present
counting argument supplies its own elementary steps. Their titles in
the bibliography identify the original sources. No knowledge of their
algebraic terminology is required to read the proof above.

\paragraph{Sources for earlier investigations of weaker bounds.}
The earlier research drafts consulted or cited the completion and
decomposition results of Barber, K\"uhn, Lo, Osthus and
Taylor~\cite{BKLOT}, Bowditch and Dukes~\cite{BD},
Yu and Feng~\cite{YF}, He and Cheng~\cite{HC}, and
Allsop, Bowtell, Lesgourgues and Petrova~\cite{ABLP}.
These references record a previous approach and their theorems are not
inputs to the counting proof in this document.


\begin{thebibliography}{99}
\bibitem{Nelder77}
J. Nelder,
\emph{Critical sets in Latin squares},
CSIRO Division of Mathematics and Statistics Newsletter 38 (1977), p.~4.
Original note; bibliographic details checked through~\cite{BeanThesis}.

\bibitem{BeanThesis}
R. W. Bean,
\emph{Critical Sets in Latin Squares and Associated Structures},
Ph.D. thesis, University of Queensland, 2001, Section~2.6.
\href{https://arxiv.org/abs/2606.12996}{arXiv:2606.12996} (posted 2026).

\bibitem{GHM}
M. Ghandehari, H. Hatami and E. S. Mahmoodian,
\emph{On the size of the minimum critical set of a Latin square},
Discrete Mathematics 293 (2005), 121--127.
\href{https://doi.org/10.1016/j.disc.2004.08.026}{doi:10.1016/j.disc.2004.08.026}.

\bibitem{BeanBCC}
P. J. Cameron,
\emph{Research problems from the 19th British Combinatorial Conference},
Discrete Mathematics 293 (1--3) (2005), 313--320.
Problem~2 (BCC19.2), proposed by R. Bean, especially part~(2).
\href{https://doi.org/10.1016/j.disc.2005.02.006}{doi:10.1016/j.disc.2005.02.006}.
\href{https://webspace.maths.qmul.ac.uk/p.j.cameron/bcc/bccp19.pdf}{Problem list}.

\bibitem{CvR}
D. Curran and G. H. J. van Rees,
\emph{Critical sets in Latin squares},
Congressus Numerantium 22 (1979), 165--168,
Proceedings of the Eighth Manitoba Conference on Numerical Mathematics
and Computing (1978).

\bibitem{SvR82}
D. R. Stinson and G. H. J. van Rees,
\emph{Some large critical sets},
Congressus Numerantium 34 (1982), 441--456.

\bibitem{Smetaniuk}
B. Smetaniuk,
\emph{On the minimal critical set of a Latin square},
Utilitas Mathematica 16 (1979), 97--100.

\bibitem{CDS}
J. Cooper, D. Donovan and J. Seberry,
\emph{Latin squares and critical sets of minimal size},
Australasian Journal of Combinatorics 4 (1991), 113--120.
\href{https://ajc.maths.uq.edu.au/pdf/4/ocr-ajc-v4-p113.pdf}{Journal text}.

\bibitem{BvR99}
J. A. Bate and G. H. J. van Rees,
\emph{The size of the smallest strong critical set in a Latin square},
Ars Combinatoria 53 (1999), 73--83.

\bibitem{BvR03}
J. A. Bate and G. H. J. van Rees,
\emph{Minimal and near-minimal critical sets in back-circulant Latin squares},
Australasian Journal of Combinatorics 27 (2003), 47--61.
\href{https://ajc.maths.uq.edu.au/pdf/27/ajc_v27_p047.pdf}{Journal text}.

\bibitem{Cavenagh}
N. J. Cavenagh,
\emph{Latin trade algorithms and the smallest critical set in a Latin square},
Journal of Automata, Languages and Combinatorics 8 (2003), 567--578.

\bibitem{BeanTrades}
R. Bean,
\emph{Latin trades on three or four rows},
Discrete Mathematics 306 (2006), 3028--3041.
\href{https://richardbean.id.au/papers/four.pdf}{Author-hosted manuscript}.

\bibitem{HQ}
H. Hatami and Y. Qian,
\emph{Teaching dimension, VC dimension, and critical sets in Latin squares},
Journal of Combinatorics 9 (2018), 9--20.
\href{https://arxiv.org/abs/1606.00032}{arXiv:1606.00032}.

\bibitem{Bala90}
K. Balasubramanian,
\emph{On transversals in Latin squares},
Linear Algebra and its Applications 131 (1990), 125--129.
\href{https://doi.org/10.1016/0024-3795(90)90378-P}{doi:10.1016/0024-3795(90)90378-P}.

\bibitem{BrunsVetter}
W. Bruns and U. Vetter,
\emph{Determinantal Rings}, Lecture Notes in Mathematics 1327,
Springer, 1988, Theorem~2.10 and Remark~2.12.
\href{https://www.home.uni-osnabrueck.de/wbruns/brunsw/detrings.pdf}{Author-hosted text}.

\bibitem{MSV}
L. E. Miller, A. K. Singh and M. Varbaro,
\emph{The $F$-pure threshold of a determinantal ideal},
Bulletin of the Brazilian Mathematical Society, New Series 45 (2014),
767--775.
\href{https://arxiv.org/abs/1210.6729}{arXiv:1210.6729}.

\bibitem{HH}
M. Hochster and C. Huneke,
\emph{Comparison of symbolic and ordinary powers of ideals},
Inventiones Mathematicae 147 (2002), 349--369.
\href{https://arxiv.org/abs/math/0211174}{arXiv:math/0211174}.

\bibitem{BKLOT}
B. Barber, D. K\"uhn, A. Lo, D. Osthus and A. Taylor,
\emph{Clique decompositions of multipartite graphs and completion of Latin squares},
Journal of Combinatorial Theory, Series A 151 (2017), 146--201.
\href{https://arxiv.org/abs/1603.01043}{arXiv:1603.01043}.

\bibitem{BD}
F. C. Bowditch and P. J. Dukes,
\emph{Fractional triangle decompositions of dense 3-partite graphs},
Journal of Combinatorics 10 (2019), 255--282.
\href{https://arxiv.org/abs/1510.08998}{arXiv:1510.08998}.

\bibitem{YF}
S. Yu and T. Feng,
\emph{On the completion of $\varepsilon$-dense partial Latin squares},
2026 preprint.
\href{https://arxiv.org/abs/2606.31143}{arXiv:2606.31143}.

\bibitem{HC}
H. He and G. Cheng,
\emph{Fractional triangle decompositions at partite minimum degree $4n/5$},
2026 preprint.
\href{https://arxiv.org/abs/2609.31758}{arXiv:2609.31758}.

\bibitem{ABLP}
J. Allsop, C. Bowtell, T. Lesgourgues and K. Petrova,
\emph{Improved bounds on completion of partial Latin squares},
2026 preprint.
\href{https://arxiv.org/abs/2609.37877}{arXiv:2609.37877}.
\end{thebibliography}
\end{document}